%% file: paper.tex
\documentclass[journal,twoside,web]{ieeecolor}

\def\logoname{LOGO-lcsys-web} 
\definecolor{subsectioncolor}{rgb}{0.067,0.627,0.859}
\def\journalname{IEEE Control Systems Letters}

\usepackage{cite}
\usepackage{epsfig}
\usepackage{tikz}           
\usetikzlibrary{shapes.geometric, arrows.meta, positioning, fit, backgrounds, calc, math, chains}    
\usepackage{pgfplots, pgfmath}
\usepgfplotslibrary{statistics, groupplots, external}
\pgfplotsset{compat=1.18}
\usepackage{graphicx}                               
\usepackage{comment}                                
\usepackage{lipsum}                                 
\usepackage{amsfonts}                               
\usepackage[utf8]{inputenc}                         
\usepackage{amsmath,amssymb,amsfonts,mathtools}     
\usepackage{url}                                    
\usepackage{nameref}                                
\usepackage{cases}                                  
\usepackage{xfrac}                                  
\usepackage{algorithm}                              
\usepackage{algpseudocodex}                         
\usepackage{float}                                  
\usepackage{times}                                  
\usepackage{colortbl}                               
\usepackage{soul}                                   
\usepackage[dvipsnames, svgnames]{xcolor}

\def\BibTeX{{\rm B\kern-.05em{\sc i\kern-.025em b}\kern-.08em
    T\kern-.1667em\lower.7ex\hbox{E}\kern-.125emX}}
\colorlet{cfn_color}{BrickRed}
\colorlet{dgn_color}{RoyalBlue}
\colorlet{input_color}{BlueViolet}
\colorlet{state_color}{RedOrange}

\newtheorem{theorem}{Theorem}
\newtheorem{proposition}{Proposition}

\newtheorem{definition}{Definition}
\newtheorem{assumption}{Assumption}

\newtheorem{remark}{Remark}

\newcommand{\norm}[2]{ \left\Vert {#1} \right\Vert_{#2}}

\newcommand{\infnorm}[1]{ \left\Vert {#1} \right\Vert}
\newcommand{\abs}[1]{ \left| {#1} \right|}

\newcommand{\defeq}{:=} 

\newcommand{\paragraphfont}[1]{\textsc{{#1}.}}
\newcommand{\figurename}{Fig.}

\newcommand{\minibatch}{\mathrm{b}}

\newcommand{\hiddenunits}{h}
\newcommand{\layers}{\mathcal{L}}
\newcommand{\sequence}{v}
\newcommand{\totsequence}{V}
\newcommand{\sequencep}{(\sequence)}
\newcommand{\parameters}{\theta}

\newcommand{\fit}{\textrm{Fit}}

\newcommand{\dISS}{$\delta$ISS}
\newcommand{\tanhfun}{\phi}
\newcommand{\realset}{\mathbb{R}}
\newcommand{\naturalset}{\mathbb{N}} 
\newcommand{\posintegerset}{\mathbb{Z}^+} 

\newlength{\customdisplayspace}
\newlength{\customjotspace}
\newenvironment{tequation}
 {\begingroup
 \let\displaystyle\textstyle
 \setlength{\abovedisplayskip}{\customdisplayspace}
 \setlength{\belowdisplayskip}{\customdisplayspace}
 \setlength{\abovedisplayshortskip}{\customdisplayspace}
 \setlength{\belowdisplayshortskip}{\customdisplayspace}
 \equation}
 {\endequation
 \endgroup}
 
\newenvironment{tequation*}
 {\begingroup
 \let\displaystyle\textstyle
 \setlength{\abovedisplayskip}{\customdisplayspace}
 \setlength{\belowdisplayskip}{\customdisplayspace}
 \setlength{\abovedisplayshortskip}{\customdisplayspace}
 \setlength{\belowdisplayshortskip}{\customdisplayspace}
 \csname equation*\endcsname}
 {\endequation
 \endgroup}

\newenvironment{talign}
 {\begingroup
 \setlength{\abovedisplayskip}{\customdisplayspace}
 \setlength{\belowdisplayskip}{\customdisplayspace}
 \setlength{\abovedisplayshortskip}{\customdisplayspace}
 \setlength{\belowdisplayshortskip}{\customdisplayspace}
 \setlength{\jot}{\customjotspace}
 \align}
 {\endalign
 \endgroup}
 
\newenvironment{talign*}
 {\begingroup
 \setlength{\abovedisplayskip}{\customdisplayspace}
 \setlength{\belowdisplayskip}{\customdisplayspace}
 \setlength{\abovedisplayshortskip}{\customdisplayspace}
 \setlength{\belowdisplayshortskip}{\customdisplayspace}
  \setlength{\jot}{\customjotspace}
 \csname align*\endcsname}
 {\endalign
 \endgroup}

\newenvironment{taligned}
 {\begingroup
  \let\displaystyle\textstyle
  \setlength{\jot}{\customjotspace}
  \aligned}
 {\endaligned
  \endgroup}

\title{Chaos-Free Networks are Stable Recurrent Neural Networks}

\author{Stefano De Carli$^{a}$, Davide Previtali$^{a}$, Mirko Mazzoleni$^{a}$ and Fabio Previdi$^{a}$%
\thanks{
$^{a}$
Department of Management, Information and Production Engineering, University of Bergamo, 24044 Dalmine, Bergamo, Italy. 
Contact: \texttt{\{name\}.\{surname\}@unibg.it} (e.g., \texttt{stefano.decarli@unibg.it}).
}
}

\begin{document}

\maketitle
\thispagestyle{empty}
\pagestyle{empty}

\begin{abstract}
Gated Recurrent Neural Networks (RNNs) are widely used for nonlinear system identification due to their high accuracy, although they often exhibit complex, chaotic dynamics that are difficult to analyze. 
This paper investigates the system-theoretic properties of the Chaos-Free Network (CFN), an architecture originally proposed to eliminate the chaotic
behavior found in standard gated RNNs. 
First, we formally prove that the CFN satisfies Input-to-State Stability (ISS) by design. 
However, we demonstrate that the CFN architecture does not intrinsically guarantee Incremental ISS (\dISS{}), as ensuring this property relies on specific parametric constraints. 
To address this, we introduce the Decoupled-Gate Network (DGN), a novel structural variant of the CFN that removes internal state connections in the gating mechanisms. 
Finally, we prove that the DGN unconditionally satisfies the \dISS{} property, providing an incrementally stable architecture for identifying nonlinear dynamical systems without requiring complex network training modifications.
Numerical results confirm that the DGN maintains the modeling capabilities of standard architectures while adhering to these rigorous stability guarantees.
\end{abstract}
%
\section{Introduction}
\label{s:introduction}
\IEEEPARstart{R}{ecurrent} Neural Networks (RNNs) have become a standard tool for the identification of nonlinear dynamical systems~\cite{pillonettoDeepNetworksSystem2025}.
While gated architectures, such as Long Short-Term Memory (LSTM) and Gated Recurrent Unit (GRU) networks~\cite[Ch. 15]{murphyProbabilisticMachineLearning2025}, demonstrate superior accuracy when it comes to modeling real systems~\cite{pillonettoDeepNetworksSystem2025}, their deployment in safety-critical control frameworks is limited by their complexity.
Specifically, standard gated RNNs can exhibit complex, irregular dynamics such as non-divergent but non-vanishing oscillations or chaotic behaviors, which are heavily influenced by the networks' initial conditions~\cite{laurent_recurrent_2016}.
Further, the lack of intrinsic stability guarantees makes the derivation of theoretically sound RNN-based control strategies challenging~\cite{bonassiRecurrentNeuralNetworks2022}.
Consequently, the analysis of RNN stability has received significant attention in recent years~\cite{terzi_learning_2021,bonassi_stability_2021,damicoIncrementalInputtoStateStability2024,decarliInfinitynormbasedInputtoStateStableLong2025}, particularly regarding \emph{Input-to-State Stability} (ISS) and \emph{Incremental Input-to-State Stability} (\dISS{}), which are crucial for robust control frameworks and observer designs~\cite{bonassiRecurrentNeuralNetworks2022,schimpernaRobustOffsetFreeConstrained2024,schimpernaRobustConstrainedNonlinear2024}.
In parallel with these control-theoretic developments, the machine learning community has actively investigated the equilibrium properties of neural networks~\cite{millerStableRecurrentModels2018a, bai2019deep} to ensure well-behaved forward propagation.
In this context, the \emph{Chaos-Free Network} (CFN) was proposed in the machine learning literature~\cite{laurent_recurrent_2016} to avoid the emergence of \emph{strange} RNN state attractors, which are chaotic fractal structures typical of standard gated networks~\cite{laurent_recurrent_2016},~\cite[Ch. 9]{strogatzNonlinearDynamicsChaos2019}, ensuring instead that the CFN's state relaxes to a stable equilibrium when no input sequence is applied to it.
While originally motivated by qualitative dynamical behavior, rigorous system-theoretic properties of the CFN remain unexplored.

In this work, we analyze the CFN through the lens of control theory. 
First, we formally prove that the absence of chaos in CFNs translates to the satisfaction of the ISS property by design. 
Second, we identify that the recurrent gating structure in CFNs does not yield an unconditional \dISS{} guarantee by design.
Indeed, promoting \dISS{} requires enforcing complex parametric conditions during CFN training, analogous to those established for LSTM and GRU networks~\cite{terzi_learning_2021,bonassi_stability_2021,decarliInfinitynormbasedInputtoStateStableLong2025}. 
To address this, inspired by sequence modeling trends that favor simplified gating mechanisms for computational efficiency~\cite{fengWereRNNsAll2024}, we propose the \emph{Decoupled-Gate Network} (DGN), a structural variant of the CFN that removes the state-dependent recurrence in the gating dynamics.
While similar architectural designs have been recently explored in the machine learning community to enhance parallelization~\cite{fengWereRNNsAll2024}, we show that this structural choice provides an important system-theoretic benefit: it renders the DGN \emph{unconditionally} \dISS{} by design.
Therefore, we establish a flexible neural network architecture for system identification that eliminates the need for constrained training procedures (e.g., stability penalties~\cite{bonassiRecurrentNeuralNetworks2022}) or complex post-training verification.

The remainder of this work is structured as follows. 
Section~\ref{s:preliminaries} introduces the notation and problem statement, and recalls the definitions of ISS and \dISS{}. 
Section~\ref{s:models} details the dynamics of the CFN and introduces the proposed DGN architecture. 
Section~\ref{s:stability} provides the main theoretical contributions of the paper, deriving the stability proofs for both architectures.
Section~\ref{s:results} presents numerical experiments validating the theoretical findings and comparing the modeling performance of the CFN and DGN against stable standard baselines.
Finally, Section~\ref{s:conclusion} draws conclusions and outlines future research directions.

\section{Preliminaries and problem statement}
\label{s:preliminaries}
\subsection{Notation and preliminaries}
\label{ss:notation_preliminaries}
Let $\realset$, $\realset_{\geq0}$, $\realset_{>0}$, $\naturalset$, $\posintegerset$ denote the sets of real, non-negative real, positive real, natural numbers (including zero), and positive integers, respectively.
Given $n, m \in \posintegerset$, $\realset^n$ is the set of real column vectors of dimension $n$, while $\realset^{n \times m}$ is the set of real matrices of dimension $n \times m$.
Vectors are typeset in boldface lowercase, e.g., $\boldsymbol{v}=\left[ v_1,\dots, v_n\right]^\top \in \realset^n$.
$\boldsymbol{0}$ and $I$ denote the zero and identity matrices of appropriate dimensions.
$\abs{\mathcal{S}}$ is the cardinality of a set $\mathcal{S}$.
The Hadamard (element-wise) product is indicated by $\circ$, and $\norm{\cdot}{p}$ represents the $p$-norm; unless otherwise specified, $\infnorm{\cdot}$ denotes the infinity norm (i.e., $p=\infty$) for vectors and the induced infinity norm for matrices.
We represent scalar discrete-time signals as $s_k$, where $k \in \naturalset$ is the time step.
Accordingly, an $n$-dimensional discrete-time signal is defined as $\boldsymbol{x}_k = \left[x_{1, k}, \dots, x_{n, k}\right]^\top \in \realset^n$.
The sequence of values attained by $\boldsymbol{x}_k$ between $k = k_1$ and $k = k_2$ is denoted by $\{ \boldsymbol{x}_k \}_{k=k_1}^{k_2}$.
We use standard activation functions, i.e., the sigmoid $\sigma (x) = \left(1 + e^{-x}\right)^{-1} \in (0,1)$ and the hyperbolic tangent $\tanhfun(x) = \tanh(x) \in (-1,1)$ for $x \in \realset$.
For a scalar function $f:\realset\to\realset$, its bold counterpart $\boldsymbol{f}:\realset^n\to\realset^n$ denotes the component-wise application.
Finally, we employ standard definitions for comparison functions of classes $\mathcal{K}$, $\mathcal{K}_{\infty}$, and $\mathcal{KL}$ as detailed in~\cite{kellettCompendiumComparisonFunction2014}.

\subsection{Problem statement}
\label{ss:problem_statement}
We consider the problem of identifying a nonlinear discrete-time Multiple-Input Multiple-Output (MIMO) dynamical system from data.
The objective is to learn the unknown map governing the evolution of the system output $\boldsymbol{y}_k \in \realset^{n_y}$ based on the input $\boldsymbol{u}_k \in \realset^{n_u}$, with $n_y, n_u \in \posintegerset$.
The system is characterized by an internal state $\boldsymbol{x}_k \in \realset^{n_x}$, where $n_x \in \posintegerset$.
We assume that the input $\boldsymbol{u}_k \in \mathcal{U} \subset \realset^{n_u}$ and the state $\boldsymbol{x}_k \in \mathcal{X} \subset \realset^{n_x}$ belong to compact sets.
The initial state $\boldsymbol{x}_0 \in \mathcal{X}$ is unknown and not subject to estimation.
To identify the underlying dynamics, we use a dataset $\mathcal{D} = \{ \mathcal{D}^{\sequencep} \}_{\sequence = 1}^{\totsequence}$ comprising $\totsequence \in \posintegerset$ input-output sequences.
Each sequence is defined as $\mathcal{D}^{\sequencep} \defeq \{ (\boldsymbol{u}_k^{\sequencep}, \boldsymbol{y}_k^{\sequencep}) \}_{k=0}^{N^{\sequencep}-1}$, where $N^{\sequencep} \in \posintegerset$ represents the length of the $\sequence$-th sequence, with $\sequence \in \{1, \ldots, \totsequence\}$.
Our goal is to identify the underlying dynamics using RNN architectures, focusing on models satisfying ISS and \dISS{}.
Stability is defined for a generic system $\boldsymbol{x}_{k+1} = \boldsymbol{w}(\boldsymbol{x}_k, \boldsymbol{u}_k)$, with $\boldsymbol{w}: \realset^{n_x} \times \realset^{n_u} \to \realset^{n_x}$, assuming that $\mathcal{X}$ is forward invariant with respect to $\mathcal{U}$, i.e., $\boldsymbol{w}(\boldsymbol{x}_k, \boldsymbol{u}_k) \in \mathcal{X}$, $\forall k \in \naturalset$, for any $\boldsymbol{x}_k \in \mathcal{X}$ and $\boldsymbol{u}_k \in \mathcal{U}$.
\begin{definition}[ISS~\cite{jiang_input--state_2001, terzi_learning_2021}]
    \label{def:iss}
    A system $\boldsymbol{x}_{k+1} = \boldsymbol{w}(\boldsymbol{x}_k, \boldsymbol{u}_k)$ is ISS in $\mathcal{X}$ with respect to $\mathcal{U}$ if there exist functions $\beta \in \mathcal{KL}$ and $\gamma_u, \gamma_b \in \mathcal{K}_\infty$ such that, for any $k \in \naturalset$, any initial state $\boldsymbol{x}_0 \in \mathcal{X}$, any input sequence $\{\boldsymbol{u}_z \in \mathcal{U} \}_{z=0}^{k-1}$, and any bias vector $\boldsymbol{b} \in \realset^{n_x}$, it holds that:
    \begin{tequation}
        \label{eq:iss}
        \infnorm{\boldsymbol{x}_k} \leq  \beta(\infnorm{\boldsymbol{x}_0}, k) + \gamma_u \left(\max_{0 \leq z < k} \infnorm{\boldsymbol{u}_z} \right) + \gamma_b \left( \infnorm{\boldsymbol{b}} \right).
    \end{tequation}
\end{definition}
\begin{definition}[\dISS{}~\cite{bayerDiscretetimeIncrementalISS2013, terzi_learning_2021}]
    \label{def:diss}
    A system $\boldsymbol{x}_{k+1} = \boldsymbol{w}(\boldsymbol{x}_k, \boldsymbol{u}_k)$ is \dISS{} in $\mathcal{X}$ with respect to $\mathcal{U}$ if there exist functions $\beta_{\delta}\in\mathcal{KL}$ and $\gamma_{\delta u}\in\mathcal{K}_{\infty}$ such that, for any $k\in\naturalset$, any pair of initial states $\boldsymbol{x}_{0}^{\mathrm{a}}, \boldsymbol{x}_{0}^{\mathrm{b}} \in\mathcal{X}$, and any pair of input sequences $\{\boldsymbol{u}_{z}^{\mathrm{a}}\in\mathcal{U}\}_{z=0}^{k-1}$ and $\{\boldsymbol{u}_{z}^{\mathrm{b}}\in\mathcal{U}\}_{z=0}^{k-1}$, the corresponding state trajectories $\boldsymbol{x}_{k}^{\mathrm{a}}$ and $\boldsymbol{x}_{k}^{\mathrm{b}}$ satisfy:
    \begin{tequation}
        \label{eq:diss}
        \infnorm{\boldsymbol{x}_{k}^{\mathrm{a}} \! - \! \boldsymbol{x}_{k}^{\mathrm{b}}} \! \leq \! \beta_{\delta} \! \left(\infnorm{ \boldsymbol{x}_{0}^{\mathrm{a}} \! - \! \boldsymbol{x}_{0}^{\mathrm{b}}},k\right) + \gamma_{\delta u} \! \left(\max_{0\leq z<k}\infnorm{ \boldsymbol{u}_{z}^{\mathrm{a}} \! -  \! \boldsymbol{u}_{z}^{\mathrm{b}}}\right) \!\!.
    \end{tequation}
\end{definition}
We recall that \dISS{} is a stronger property than ISS as the class of \dISS{} systems constitutes a strict subset of the class of ISS systems~\cite[Proposition 1]{bayerDiscretetimeIncrementalISS2013}.
\begin{remark}
    \label{rem:inf_norm}
    In this work, ISS and \dISS{} are defined using the infinity norm for convenience, but the properties hold for any other vector norm~\cite{jiang_input--state_2001}.
\end{remark}

\section{Chaos-free and decoupled-gate networks}
\label{s:models}
While standard LSTM and GRU networks are widely used and offer great modeling ability~\cite{pillonettoDeepNetworksSystem2025},~\cite[Ch. 15]{murphyProbabilisticMachineLearning2025}, they often exhibit chaotic behaviors~\cite{laurent_recurrent_2016} and lack theoretical stability guarantees by design.
To address this, we analyze the chaos-free network architecture~\cite{laurent_recurrent_2016}, which mitigates these issues, and propose a novel CFN variant, the decoupled-gate network, designed to guarantee unconditional \dISS{}.
\subsection{The Chaos-Free Network (CFN)}
\label{ss:cfn_standard}
The CFN was introduced in~\cite{laurent_recurrent_2016} for general machine learning tasks.
CFNs are designed to mitigate the chaotic behavior exhibited by standard gated RNNs (LSTM and GRU networks) while still managing the vanishing gradient problem.
The network consists of $L \in \posintegerset$ stacked layers.
Each layer $l \in \layers = \{1, \ldots, L\}$ is a nonlinear discrete-time MIMO dynamical system in state-space form with hidden state $\boldsymbol{h}_k^{(l)} \in \realset^{n_{\hiddenunits{}}^{(l)}}$.
Here, $(l)$ denotes the $l$-th layer, and $n_{\hiddenunits{}}^{(l)} \in \posintegerset$ is its number of hidden units.
The hidden state evolution of each $l$-th layer is governed by a forget gate $\boldsymbol{f}_k^{(l)}  \in \realset^{n_{\hiddenunits{}}^{(l)}}$, an input gate $\boldsymbol{i}_k^{(l)} \in  \realset^{n_{\hiddenunits{}}^{(l)}}$, and a candidate hidden state $\boldsymbol{\tilde{h}}_k^{(l)} \in  \realset^{n_{\hiddenunits{}}^{(l)}}$:
\begin{subequations}
    \label{eq:CFN_general_dynamics}
    \begin{talign}
        \label{eq:CFN_layer_state_update}
        \boldsymbol{h}_{k+1}^{(l)} &\defeq \boldsymbol{f}_{k}^{(l)} \circ \boldsymbol{\tanhfun}\left(\boldsymbol{h}_{k}^{(l)}\right) + \boldsymbol{i}_{k}^{(l)} \circ \boldsymbol{\tilde{h}}_{k}^{(l)}, \\
        \label{eq:gates_forget}
        \boldsymbol{f}_k^{(l)} &\defeq \boldsymbol{\sigma}\left( W_f^{(l)} {\boldsymbol{\tilde{u}}}_k^{(l)}  +  R_f^{(l)} \boldsymbol{h}_k^{(l)}  +  \boldsymbol{b}_f^{(l)} \right), \\
        \label{eq:gates_input}
        \boldsymbol{i}_k^{(l)} &\defeq \boldsymbol{\sigma}\left( W_i^{(l)} {\boldsymbol{\tilde{u}}}_k^{(l)}  +  R_i^{(l)} \boldsymbol{h}_k^{(l)}  +  \boldsymbol{b}_i^{(l)} \right), \\
        \label{eq:gates_candidate}
        \boldsymbol{\tilde{h}}_k^{(l)} &\defeq \boldsymbol{\tanhfun}\left( W_{\tilde{h}}^{(l)} {\boldsymbol{\tilde{u}}}_k^{(l)}  + \boldsymbol{b}_{\tilde{h}}^{(l)} \right), \\
        \label{eq:layer_input_def}
        \boldsymbol{\tilde{u}}_k^{(l)} &\defeq
        \begin{cases}
            \boldsymbol{u}_k & \text{if } l = 1, \quad n_{\tilde{u}}^{(l)} = n_u, \\
            \boldsymbol{h}_{k+1}^{(l-1)} & \text{if } l \in \layers \setminus \{1\}, \quad n_{\tilde{u}}^{(l)} = n_{\hiddenunits}^{(l-1)},
        \end{cases}
    \end{talign}
\end{subequations}
where $\boldsymbol{\tilde{u}}_k^{(l)} \in \realset^{n_{\tilde{u}}^{(l)}}$ is the layer input, $W_f^{(l)}, W_i^{(l)}, W_{\tilde{h}}^{(l)} \in \realset^{n^{(l)}_\hiddenunits{} \times n_{\tilde{u}}^{(l)}}$ are the input weights, $R_f^{(l)}, R_i^{(l)} \in \realset^{n^{(l)}_\hiddenunits{} \times n^{(l)}_\hiddenunits{}}$ are the recurrent weights, and $\boldsymbol{b}_f^{(l)}, \boldsymbol{b}_i^{(l)}, \boldsymbol{b}_{\tilde{h}}^{(l)} \in \realset^{n_{\hiddenunits{}}^{(l)}}$ are the bias vectors.
These single-layer CFN dynamics are shown in \figurename{}~\ref{fig:CFN-DGN}.
A key aspect of CFNs is that, due to the ranges of the activation functions $\sigma(\cdot) \in (0,1)$ and $\tanhfun(\cdot) \in (-1,1)$, every layer admits a compact forward invariant set.
\begin{proposition}[Forward invariant set of $\boldsymbol{h}_k^{(l)}$]
    \label{pr:invariant_state_layer}
    For any $l \in \layers$, the set $\mathcal{H}_{\mathrm{inv}}^{(l)} \defeq \left[-2, 2\right]^{n_{\hiddenunits}^{(l)}}$ is a forward invariant compact set for the layer dynamics in~\eqref{eq:CFN_general_dynamics}.
    Thus, for any initial condition $\boldsymbol{h}_{0}^{(l)} \in \mathcal{H}_{\mathrm{inv}}^{(l)}$, the hidden state satisfies $\boldsymbol{h}_{k}^{(l)} \in \mathcal{H}_{\mathrm{inv}}^{(l)}$ for all $k \in \naturalset$ and any input $\boldsymbol{\tilde{u}}_k^{(l)} \in \realset^{n_{\tilde{u}}^{(l)}}$.
\end{proposition}
The proof is reported in Appendix~\ref{s:proofs}.
In practice, $\boldsymbol{h}_{0}^{(l)}$ is typically initialized as the zero vector~\cite[Ch. 10]{goodfellowDeepLearning2016}, i.e., $\boldsymbol{h}_{0}^{(l)} = \boldsymbol{0} \in \mathcal{H}_{\mathrm{inv}}^{(l)}, \forall l \in \layers$ .
By Proposition~\ref{pr:invariant_state_layer}, the complete network state $\boldsymbol{h}_k = [{\boldsymbol{h}_{k}^{(1)}}^\top, \dots, {\boldsymbol{h}_{k}^{(L)}}^\top]^\top$ evolves within $\mathcal{H}_{\mathrm{inv}} = \bigtimes_{l=1}^L \mathcal{H}_{\mathrm{inv}}^{(l)} \subset \realset^{n_\hiddenunits{}}$, where $n_\hiddenunits{} = \sum_{l=1}^L n_\hiddenunits^{(l)}$ is the total number of hidden units.
The network output $\boldsymbol{y}_k$ follows from a fully connected mapping of the final layer hidden state:
\begin{tequation}
    \label{eq:CFN}
    \left\{\begin{taligned}
    \boldsymbol{h}_{k+1}^{(l)} &\defeq \boldsymbol{f}_{k}^{(l)} \circ \boldsymbol{\tanhfun}\left(\boldsymbol{h}_{k}^{(l)}\right) + \boldsymbol{i}_{k}^{(l)} \circ \boldsymbol{\tilde{h}}_{k}^{(l)}, \   \forall l \in \layers, \\
    \boldsymbol{y}_k &\defeq W_y \boldsymbol{h}_{k+1}^{(L)} + \boldsymbol{b}_y,
    \end{taligned}\right.
\end{tequation}
where $W_y \in \realset^{n_y \times n_{\hiddenunits}^{(L)}}$ and $\boldsymbol{b}_y \in \realset^{n_y}$.
Finally, the CFN relies on the following set of tunable parameters $\parameters$:
\begin{tequation}
\begin{taligned}
    \label{eq:CFN_parameters}
    \parameters &\defeq{} \{W_y, \boldsymbol{b}_y\}  \cup  \{\parameters^{(l)}: l \in \layers \}, \\
    %
    \parameters^{(l)} &\defeq{} \{W_f^{(l)}, W_{i}^{(l)}, W_{\tilde{h}}^{(l)}, R_f^{(l)}, R_{i}^{(l)}, \boldsymbol{b}_f^{(l)}, \boldsymbol{b}_i^{(l)}, \boldsymbol{b}_{\tilde{h}}^{(l)}\}.
\end{taligned}
\end{tequation}
A distinctive feature of CFNs, compared to common gated RNNs, is that the candidate state $\boldsymbol{\tilde{h}}_k^{(l)}$, $l \in \layers$, depends \emph{solely} on the current layer input $\boldsymbol{\tilde{u}}_k^{(l)}$ and not on the hidden state at the previous time step.
However, the CFN retains full recurrent connectivity due to the presence of the weight matrices $R_f^{(l)}, R_i^{(l)}$.
This coupling between the hidden state $\boldsymbol{h}_k^{(l)}$ and the gates introduces a state-dependent modulation that, while enabling rich dynamics, makes it challenging to provide a priori strong stability properties (such as \dISS{}) independently of the specific parameter values, see Section~\ref{s:stability}.
\subsection{The proposed Decoupled-Gate Network (DGN)}
\label{ss:dgn_proposed}
To guarantee unconditional \dISS{} without parametric constraints, we introduce the decoupled-gate network.
Inspired by parallelizable RNN architectures~\cite{fengWereRNNsAll2024}, from which it inherits training efficiency benefits, the DGN structurally simplifies the CFN by eliminating the recurrent dependencies of the gates on the hidden state, setting $R_f^{(l)} = R_i^{(l)} = \boldsymbol{0}$ for all $l \in \layers$.
Consequently, the gates in~\eqref{eq:gates_forget} and~\eqref{eq:gates_input} of a DGN layer depend exclusively on the layer input $\boldsymbol{\tilde{u}}_k^{(l)}$, removing state-dependent modulations that complicate stability analysis.
Although this restricts single-layer dynamics, complex temporal processing is recovered through network depth ($L > 1$)~\cite{fengWereRNNsAll2024}.
Since deeper layers are driven by previous states ($\boldsymbol{\tilde{u}}_k^{(l)} = \boldsymbol{h}_{k+1}^{(l-1)}$), their gating mechanisms remain effectively state-dependent.
This structural choice preserves model expressivity while rendering the DGN \dISS{} \emph{by design}, as proven in Section~\ref{s:stability}.
\begin{figure}
    \centering
    \input{figures/cfn-dgn}
    \caption{Single-layer dynamics for the CFN and the DGN.
    The dashed paths represent the recurrent gate connections, which are absent in DGNs.
    Orange and purple represent the hidden state and input paths, respectively.
    }
    \label{fig:CFN-DGN}
\end{figure}
%
\section{Stability properties}
\label{s:stability}
In this section, we analyze the stability of the CFN and DGN.
All proofs are deferred to Appendix~\ref{s:proofs}.
We adopt the following assumption on input boundedness.
\begin{assumption}[Input boundedness]
    \label{as:input_boundedness}
    The inputs $\boldsymbol{u}_{k}$ are bounded such that $\boldsymbol{u}_{k} \in \mathcal{U} \defeq [-1, 1]^{n_u}$ for all $k \in \naturalset$.
\end{assumption}
This is a customary assumption for neural networks, satisfied by input data normalization.
Furthermore, considering the layer input $\boldsymbol{\tilde{u}}_k^{(l)}$ in~\eqref{eq:layer_input_def}, it follows from Assumption~\ref{as:input_boundedness} and Proposition~\ref{pr:invariant_state_layer} that $\boldsymbol{\tilde{u}}_k^{(l)} \in \tilde{\mathcal{U}}^{(l)} \defeq [-2, 2]^{n_{\tilde{u}}^{(l)}}$ for all $l \in \layers$.
Indeed, $\boldsymbol{\tilde{u}}_k^{(1)} = \boldsymbol{u}_k \in \mathcal{U} \subset \tilde{\mathcal{U}}^{(1)}$, while for $l \in \layers \setminus \{ 1\}$, $\boldsymbol{\tilde{u}}_k^{(l)} = \boldsymbol{h}_{k+1}^{(l-1)} \in \mathcal{H}_{\mathrm{inv}}^{(l-1)} = \tilde{\mathcal{U}}^{(l)}$.
\subsection{Stability of CFN architectures}
\label{ss:stability_cfn}
Based on the above assumptions, we derive sufficient parametric conditions for ISS and \dISS{} of CFN layers.
\begin{theorem}[ISS of a CFN layer]
    \label{th:CFN_layer_ISS}
    Under Assumption~\ref{as:input_boundedness}, the $l$-th layer, $l \in \layers$, of the CFN in~\eqref{eq:CFN_general_dynamics} is ISS in $\mathcal{H}_{\mathrm{inv}}^{(l)}$ with respect to $\tilde{\mathcal{U}}^{(l)}$.
\end{theorem}
\begin{theorem}[\dISS{} of a CFN layer]
    \label{th:CFN_layer_dISS}
    Under Assumption~\ref{as:input_boundedness}, the $l$-th layer, $l \in \layers$, of the CFN in~\eqref{eq:CFN_general_dynamics} is \dISS{} in $\mathcal{H}_{\mathrm{inv}}^{(l)}$ with respect to $\tilde{\mathcal{U}}^{(l)}$ if:
    \begin{tequation}
        \label{eq:CFN_layer_state_update_dISS_condition}
        \rho^{(l)} \defeq \bar{\sigma}_{f}^{(l)} + \frac{1}{4} \infnorm{R_{f}^{(l)}} + \frac{1}{4} \infnorm{R_{i}^{(l)}} \bar{\tanhfun}_{\tilde{h}}^{(l)} < 1,
    \end{tequation}
    where:
    \begin{subequations}
        \begin{talign}
        \label{eq:bar_sigma_def}
        \bar{\sigma}_{f}^{(l)} &\defeq \sigma\left( \infnorm{ \begin{bmatrix} 2W_{f}^{(l)} & 2R_{f}^{(l)} & \boldsymbol{b}_{f}^{(l)} \end{bmatrix} } \right), \\
        \label{eq:bar_phi_def}
        \bar{\tanhfun}_{\tilde{h}}^{(l)} &\defeq \tanhfun \left( \infnorm{ \begin{bmatrix} 2W_{\tilde{h}}^{(l)} & \boldsymbol{b}_{\tilde{h}}^{(l)} \end{bmatrix} } \right).
         \end{talign}
    \end{subequations}
\end{theorem}
We now generalize the results to deep networks ($L > 1$).
\begin{theorem}[ISS of a CFN]
    \label{th:network_ISS}
    Under Assumption~\ref{as:input_boundedness}, the $L$-layer CFN in~\eqref{eq:CFN} is ISS in $\mathcal{H}_{\mathrm{inv}}$ with respect to $\mathcal{U}$.
\end{theorem}
\begin{theorem}[\dISS{} of a CFN]
    \label{th:network_dISS}
    Under Assumption~\ref{as:input_boundedness}, the $L$-layer CFN in~\eqref{eq:CFN} is \dISS{} in $\mathcal{H}_{\mathrm{inv}}$ with respect to $\mathcal{U}$ if Condition~\eqref{eq:CFN_layer_state_update_dISS_condition} holds for every layer $l \in \layers$.
\end{theorem}
\subsection{Extension to DGN architectures}
\label{ss:stability_DGN}
The DGN architecture, being a structural simplification of the CFN, inherits similar but stronger stability properties.
\begin{theorem}[\dISS{} of a DGN]
    \label{th:DGN_stability}
    Under Assumption~\ref{as:input_boundedness}, the $L$-layer DGN, defined by~\eqref{eq:CFN} with $R_f^{(l)} = R_i^{(l)} = \boldsymbol{0}$ for all $l \in \layers$, is \dISS{} in $\mathcal{H}_{\mathrm{inv}}$ with respect to $\mathcal{U}$.
\end{theorem}

\section{Numerical results}
\label{s:results}
In this section, we evaluate the modeling capabilities of the CFN and the proposed DGN architectures for system identification tasks, a domain \emph{unexplored} in prior works, since CFNs were validated only for word-level language modeling~\cite{laurent_recurrent_2016}.
We assess their performance on two standard benchmarks for RNN stability: the pH Neutralization Reactor~\cite{terzi_learning_2021} and the Four-Tank system~\cite{bonassi_stability_2021}.
\subsection{Training configuration}
\label{ss:training_config}
Recall the dataset $\mathcal{D}$ introduced in Section~\ref{ss:problem_statement}.
For each benchmark, we partition the set of indices $\{1, \dots, \totsequence\}$ into three disjoint subsets $\mathcal{I}_{\mathrm{tr}}$ (training), $\mathcal{I}_{\mathrm{val}}$ (validation), and $\mathcal{I}_{\mathrm{tst}}$ (test), defining the corresponding datasets $\mathcal{D}_{\mathrm{tr}}$, $\mathcal{D}_{\mathrm{val}}$, and $\mathcal{D}_{\mathrm{tst}}$.
To satisfy Assumption~\ref{as:input_boundedness}, inputs and outputs are linearly rescaled within the range $[-1, 1]$ using parameters computed exclusively from $\mathcal{D}_{\mathrm{tr}}$.
We follow the standard training procedure for RNNs~\cite[Ch. 10]{goodfellowDeepLearning2016}, also shared with baseline comparisons~\cite{terzi_learning_2021,bonassi_stability_2021}.
The parameters $\parameters$ in~\eqref{eq:CFN_parameters} are estimated by minimizing the Mean Squared Error (MSE) on $\mathcal{D}_{\mathrm{tr}}$:
\begin{tequation}
    \label{eq:mse}
    J(\parameters) \defeq \frac{1}{\abs{\mathcal{I}_{\mathrm{tr}}}} \sum_{\sequence \in \mathcal{I}_{\mathrm{tr}}}  \frac{1}{N^{\sequencep} \! - \! N_{\mathrm{w}}} \sum_{k=N_{\mathrm{w}}}^{N^{\sequencep}-1} \norm{\boldsymbol{y}_k^{\sequencep} - \boldsymbol{\hat{y}}_k^{\sequencep}(\parameters)}{2}^2,
\end{tequation}
where $\boldsymbol{\hat{y}}_k(\parameters)$ is the output predicted by the network with parameters $\parameters$ for the sequence $\sequence \in \mathcal{I}_{\mathrm{tr}}$, given the input $\boldsymbol{u}_k^{\sequencep}$ at time $k$.
We apply a washout period of $N_{\mathrm{w}}=25$ steps to the training and validation sequences, ensuring that the prediction error during the initial transient is excluded from the loss calculation~\cite[Ch. 10]{goodfellowDeepLearning2016}.
The models are trained for $2000$ epochs using the Adam optimizer with a base learning rate of $0.001$, decaying by $0.9$ every $200$ epochs to facilitate convergence.
The training data is processed in mini-batches of size $\minibatch \in \posintegerset$, which is adjusted according to the specific dataset dimensions.
To prevent overfitting, we apply a dropout rate of $0.05$.
We adopt an early stopping strategy, selecting the parameters $\parameters^*$ that yield the minimal MSE on $\mathcal{D}_{\mathrm{val}}$, which is evaluated every epoch.
Crucially, for the DGN architecture, \emph{no additional constraints are imposed during training}, as the model is \dISS{} by design (Theorem~\ref{th:DGN_stability}).
In contrast, the stable baselines (LSTM and GRU) used for comparison require minimizing a loss function augmented with a stability penalty term to promote \dISS{} compliance~\cite{terzi_learning_2021, bonassi_stability_2021}.
%
%
Lastly, we train the CFN \emph{without} auxiliary stability penalties.
Notably, the resulting CFN models consistently violated the sufficient \dISS{} condition in~\eqref{eq:CFN_layer_state_update_dISS_condition}, highlighting that \dISS{} compliance for this architecture remains parameter-dependent rather than structurally guaranteed a priori.
\subsection{pH Reactor}
\label{ss:ph_results}
We first utilize the pH neutralization benchmark~\cite{terzi_learning_2021} to compare the CFN and the DGN against the \dISS{}-LSTM network baseline established in the literature~\cite{terzi_learning_2021}, which was trained using a stability-promoting augmented loss.
The system describes the nonlinear control of the pH level ($y$) via an alkaline flow ($u$) under acid disturbances, representing a single-input single-output system.
$\mathcal{D}$ contains $\totsequence=20$ sequences generated by a Multilevel Pseudo-Random Signal (MPRS) with a sampling time of $10$ s.
$\mathcal{D}$ is partitioned into $\abs{\mathcal{D}_{\mathrm{tr}}} \! = \!15$, $\abs{\mathcal{D}_{\mathrm{val}}} \! = \!4$ and $\abs{\mathcal{D}_{\mathrm{tst}}} \!= \!1$ sequences, following the original experimental setup in~\cite{terzi_learning_2021}.

We train single-layer CFN and DGN models with $n_\hiddenunits^{(1)} = 7$ hidden units, matching the \dISS{}-LSTM model in~\cite{terzi_learning_2021}, and use a mini-batch size of $\minibatch = 5$.
Performance is evaluated for a generic sequence $\sequence \in \mathcal{I}_{\mathrm{tst}}$ using the $\fit{}^{\sequencep}$ metric:
\begin{equation*}
    \fit{}^{\sequencep}\left(\parameters^*\right) \defeq 100 \left( 1 - \frac{\sqrt{\sum_{k = 0}^{N^{\sequencep} - 1}\left( y_{k}^{\sequencep} - \hat{y}_{k}^{\sequencep} \left(\parameters^*\right) \right)^2}}{\sqrt{\sum_{k = 0}^{N^{\sequencep} - 1}\left( y_{k}^{\sequencep} - \bar{y}^{\sequencep} \right)^2}} \right),
\end{equation*}
where $\bar{y}^{\sequencep}$ is the sample mean of $y$ along sequence $\sequence$.
The results, summarized in terms of $\fit{}$ on $\mathcal{D}_{\mathrm{tst}}$, demonstrate the competitive performance of the proposed architectures, as visually confirmed by the time-domain comparison in~\figurename~\ref{fig:realizations}.
The CFN achieves a $\fit{}$ of $96.1\%$, closely matching the $96.5\%$ of the \dISS{}-LSTM benchmark in~\cite{terzi_learning_2021}. 
Notably, the DGN reaches $95.6\%$, practically matching the CFN performance. 
Therefore, the DGN maintains high accuracy while providing rigorous unconditional stability, bypassing the complex stability-promoting augmented loss functions required by standard gated RNNs~\cite{terzi_learning_2021,bonassi_stability_2021}.
\subsection{Four-Tank System}
\label{ss:four_tank_results}
We also evaluate the architectures on the Four-Tank benchmark, already used for \dISS{}-GRU network validation in~\cite{bonassi_stability_2021}.
The task is to model the water level of the lower tank ($y$) given two input pump flow rates ($\boldsymbol{u}$), constituting a multiple-input single-output problem.
A $25,\!000$-step MPRS trajectory with a sampling time of $15$ s is partitioned into $10$ blocks.
To assess the robustness of the results, we repeat the training and validation phases for all $45$ possible combinations of these blocks, selecting $8$ blocks for training and $2$ for validation in each run.
Windowing these blocks into subsequences of length $500$ (with $50\%$ overlap) practically yields $\abs{\mathcal{D}_{\mathrm{tr}}} = 77$ and $\abs{\mathcal{D}_{\mathrm{val}}} = 18$ subsequences.
Then, an independent test trajectory ($\abs{\mathcal{D}_{\mathrm{tst}}} = 1$) of $5000$ steps is generated. 
To assess performance under realistic conditions, we add white Gaussian noise to the output measurements with a signal-to-noise ratio of $100$ (i.e., $20$ dB).
We employ a 3-layer architecture with $n_{\hiddenunits}^{(l)}=7$ units per layer for both the CFN and DGN to match the \dISS{}-GRU model in~\cite{bonassi_stability_2021}, and use a mini-batch size of $\minibatch = 25$.
The models demonstrate remarkable performance on $\mathcal{D}_{\mathrm{tst}}$ for this task as well (see~\figurename~\ref{fig:realizations}, where the realization corresponding to the median $\fit{}$ performance is shown).
The aggregate results, reported as mean $\pm$ standard deviation, confirm the high robustness of the proposed models. 
Specifically, the CFN attains an average $\fit{}$ of $94.03\% \pm 0.21\%$, while the DGN follows with $92.55\% \pm 0.19\%$; notably, if the first $N_{\mathrm{w}}$ transient steps are excluded from the $\fit{}$ calculation, we obtain $96.22\% \pm 0.19\%$ for the CFN and $95.50\% \pm 0.22\%$ for the DGN.
We compare these against the \dISS{}-GRU baseline reported in~\cite{bonassi_stability_2021}, obtained with a stability-promoting augmented loss strategy.
Although this baseline achieved a $\fit{}$ of $97.0\%$, it is crucial to note that it was trained and tested with different data from the same system.

In conclusion, while structural constraints in DGNs slightly reduce expressivity, this is a favorable trade-off for obtaining unconditional \dISS{} \emph{by design}, eliminating the need for complex augmented losses during training~\cite{terzi_learning_2021, bonassi_stability_2021, decarliInfinitynormbasedInputtoStateStableLong2025}.

\begin{figure}[t]
    \centering
    \input{figures/realizations}
    \caption{Test predictions for pH Reactor (top) and Four-Tank (bottom).
    The plots compare the ground truth (dashed black line) against the CFN (thick orange line) and DGN (thin blue line) predictions.
    For the Four-Tank, median performance realizations are shown.
    }
    \label{fig:realizations}
\end{figure}
%
\section{Conclusion}
\label{s:conclusion}
In this work, we investigated the stability properties of the CFN and the proposed DGN, a novel architecture designed to ensure unconditional \dISS{} in nonlinear system identification.
We proved that the CFN inherently satisfies ISS by design, while the DGN unconditionally satisfies \dISS{} without the restrictive training constraints required by standard gated RNNs used in previous works~\cite{terzi_learning_2021,bonassi_stability_2021, decarliInfinitynormbasedInputtoStateStableLong2025}.
This contrasts with standard gated RNNs (e.g., LSTM and GRU networks), where \dISS{} must be enforced by augmenting the loss function with penalty terms during the network training process~\cite{terzi_learning_2021,bonassi_stability_2021, decarliInfinitynormbasedInputtoStateStableLong2025}.
Numerical results confirm that the DGN maintains predictive accuracies comparable to the CFN and stable gated RNN baselines while offering intrinsically strong structural stability guarantees.
Future research will focus on exploiting the intrinsic \dISS{} property of the DGN to design robust state observers and integrating these models into safety-critical control frameworks, such as model predictive control~\cite{schimpernaRobustOffsetFreeConstrained2024,schimpernaRobustConstrainedNonlinear2024}.
%
\appendices
\section{Proofs}
\label{s:proofs}
%
%

\paragraphfont{Proof of Proposition~\ref{pr:invariant_state_layer}} 
Consider the hidden state $\boldsymbol{h}_{k}^{(l)} = [h_{1,k}^{(l)}, \dots, h_{n_{\hiddenunits}^{(l)},k}^{(l)}]^\top$ for $l \in \layers$.
For any component $j \in \mathcal{J}^{(l)} \defeq \{1, \dots, n_{\hiddenunits}^{(l)}\}$, the update~\eqref{eq:CFN_layer_state_update} satisfies $|h_{j,k+1}^{(l)}| \leq |f_{j,k}^{(l)}| |\tanhfun(h_{j,k}^{(l)})| + |i_{j,k}^{(l)}| |\tilde{h}_{j,k}^{(l)}|$ by the triangle inequality.
%
%
Since $\sigma(\cdot) \in (0,1)$ and $\tanhfun(\cdot) \in (-1,1)$, it follows that $|h_{j,k+1}^{(l)}| \leq 2$.
Thus, $\forall k \in \naturalset$, if $\boldsymbol{h}_{k}^{(l)} \! \in \! \mathcal{H}_{\mathrm{inv}}^{(l)} \! \defeq \! [-2, 2]^{n_{\hiddenunits}^{(l)}}$, then $\boldsymbol{h}_{k+1}^{(l)} \! \in \! \mathcal{H}_{\mathrm{inv}}^{(l)}$ for any $\boldsymbol{\tilde{u}}_k^{(l)} \! \in \!\realset^{n_{\tilde{u}}^{(l)}}$, making $\mathcal{H}_{\mathrm{inv}}^{(l)}$ a forward invariant set.
\hfill$\blacksquare$
%

\paragraphfont{Proof of Theorem~\ref{th:CFN_layer_ISS}}
Consider a single layer $l \in \layers$ and omit the superscript $(l)$ for ease of discussion.
For any $j \in \mathcal{J}$, we can upper bound the forget gate component $f_{j,k}$ from~\eqref{eq:gates_forget} by exploiting the monotonic properties of the sigmoid function:
\begin{talign}
\label{eq:sigma_bar_deriv}
    \left|f_{j, k}\right| &\leq \max_{\tilde{\boldsymbol{u}}_k \in \tilde{\mathcal{U}}, \boldsymbol{h}_k \in \mathcal{H}_{\mathrm{inv}}} \infnorm{\boldsymbol{\sigma}\left( W_f {\boldsymbol{\tilde{u}}}_k  +  R_f \boldsymbol{h}_k  +  \boldsymbol{b}_f \right)} \nonumber \\
    &\leq \sigma\left( \infnorm{ \begin{bmatrix} 2W_{f} & 2R_{f} & \boldsymbol{b}_{f} \end{bmatrix} } \right) \defeq \bar{\sigma}_f \in (0, 1).
\end{talign}
Thus, $\infnorm{\boldsymbol{f}_k} \leq \bar{\sigma}_f$ for all $k \in \naturalset$.
Taking the norm of the state update in~\eqref{eq:CFN_layer_state_update} and applying the triangle inequality yields:
\begin{equation*}
\label{eq:iss_proof_step0}
    \infnorm{\boldsymbol{h}_{k+1}} \leq \infnorm{\boldsymbol{f}_{k}} \infnorm{\boldsymbol{\tanhfun}\left(\boldsymbol{h}_{k}\right)} + \infnorm{\boldsymbol{i}_{k}} \infnorm{\boldsymbol{\tilde{h}}_{k}}.
\end{equation*}
Using the 1-Lipschitz continuity of the hyperbolic tangent function and the bounds $\infnorm{\boldsymbol{f}_k} \leq \bar{\sigma}_f,\infnorm{\boldsymbol{i}_k} \leq 1$, we obtain:
\begin{tequation}
\label{eq:iss_proof_step1}
    \infnorm{\boldsymbol{h}_{k+1}} \leq \bar{\sigma}_{f} \infnorm{\boldsymbol{h}_{k}} + \infnorm{\boldsymbol{\tilde{h}}_k}.
\end{tequation}
From~\eqref{eq:gates_candidate}, the 1-Lipschitz continuity of the hyperbolic tangent function implies that the candidate hidden state is bounded as $\infnorm{\boldsymbol{\tilde{h}}_k} \leq \infnorm{W_{\tilde{h}}} \infnorm{\boldsymbol{\tilde{u}}_k} + \infnorm{\boldsymbol{b}_{\tilde{h}}}$.
Substituting this back into~\eqref{eq:iss_proof_step1}, we get the recursive inequality:
\begin{equation*}
\label{eq:iss_proof_step2}
     \infnorm{\boldsymbol{h}_{k+1}} \leq \bar{\sigma}_{f} \infnorm{\boldsymbol{h}_{k}} + \infnorm{W_{\tilde{h}}}\infnorm{\boldsymbol{\tilde{u}}_k} + \infnorm{\boldsymbol{b}_{\tilde{h}}}.
\end{equation*}
We propagate this inequality back to the initial state:
\begin{talign*}
    \infnorm{\boldsymbol{h}_{k}} \! &\leq \! \bar{\sigma}_{f}^k \infnorm{\boldsymbol{h}_{0}} + \sum_{z=0}^{k-1} \bar{\sigma}_f^{k-1-z} \left( \infnorm{W_{\tilde{h}}} \infnorm{\boldsymbol{\tilde{u}}_z} + \infnorm{\boldsymbol{b}_{\tilde{h}}} \right) \\
    &\leq \! \bar{\sigma}_{f}^k \! \infnorm{\boldsymbol{h}_{0}} \! + \! \left( \infnorm{W_{\tilde{h}}} \! \max_{0 \leq z < k} \! \infnorm{\boldsymbol{\tilde{u}}_z} \! + \! \infnorm{\boldsymbol{b}_{\tilde{h}}} \right) \! \sum_{z=0}^{k-1} \bar{\sigma}_f^{k-1-z}.
\end{talign*}
Since $\bar{\sigma}_{f} \in (0, 1)$, the geometric series can be bounded as $\sum_{z=0}^{k-1} \bar{\sigma}_f^{k-1-z} \leq \frac{1}{1-\bar{\sigma}_{f}}$.
We can now deduce the $\mathcal{KL}$ and the two $\mathcal{K}_{\infty}$ functions for the ISS condition in~\eqref{eq:iss}:
\begin{talign*}
    \beta\left(\infnorm{\boldsymbol{h}_{0}}, k\right) &\defeq \bar{\sigma}_{f}^k \infnorm{\boldsymbol{h}_0}, \\
    \gamma_u \left(\max_{0 \leq z < k} \infnorm{\boldsymbol{\tilde{u}}_{z}} \right) &\defeq \frac{\infnorm{W_{\tilde{h}}}}{1 - \bar{\sigma}_{f}} \max_{0 \leq z < k} \infnorm{\boldsymbol{\tilde{u}}_{z}}, \\
    \gamma_b \left(\infnorm{\boldsymbol{b}_{\tilde{h}}} \right) &\defeq \frac{1}{1 - \bar{\sigma}_{f}} \infnorm{\boldsymbol{b}_{\tilde{h}}}.
\end{talign*}
Thus, the $l$-th CFN layer, $l \in \layers$, is ISS in $\mathcal{H}_{\mathrm{inv}}^{(l)}$ with respect to $\tilde{\mathcal{U}}^{(l)}$ due to the existence of these functions.
\hfill$\blacksquare$
%

\paragraphfont{Proof of Theorem~\ref{th:CFN_layer_dISS}}
Consider a single layer $l \in \layers$ and omit the superscript $(l)$ for ease of discussion.
For any $j \in \mathcal{J}$, we first derive the bound of the candidate hidden state component $\tilde{h}_{j,k}$ from~\eqref{eq:gates_candidate}.
As with $\bar{\sigma}_{f}$ in~\eqref{eq:sigma_bar_deriv}, since the hyperbolic tangent function is monotonically increasing:
\begin{tequation}
\label{eq:tanh_bar_deriv}
    \begin{taligned}
    \left|\tilde{h}_{j,k}\right| &\leq \max_{\tilde{\boldsymbol{u}}_k \in \tilde{\mathcal{U}}} \infnorm{\boldsymbol{\tanhfun}\left(W_{\tilde{h}}\tilde{\boldsymbol{u}}_k + \boldsymbol{b}_{\tilde{h}}\right)} \\
    &\leq \tanhfun\left( \infnorm{ \begin{bmatrix} 2W_{\tilde{h}} & \boldsymbol{b}_{\tilde{h}} \end{bmatrix} } \right) \defeq{} \bar{\tanhfun}_{\tilde{h}} \in (0, 1).
    \end{taligned}
\end{tequation}
Now, consider two state trajectories $\boldsymbol{h}_k^{\mathrm{a}}, \boldsymbol{h}_k^{\mathrm{b}}$ starting from $\boldsymbol{h}_0^{\mathrm{a}}, \boldsymbol{h}_0^{\mathrm{b}} \in \mathcal{H}_{\mathrm{inv}}$ and driven by the input sequences $\{\boldsymbol{\tilde{u}}_z^{\mathrm{a}} \in \tilde{\mathcal{U}}\}_{z=0}^{k-1}$ and $\{\boldsymbol{\tilde{u}}_z^{\mathrm{b}} \in \tilde{\mathcal{U}}\}_{z=0}^{k-1}$.
Let $\boldsymbol{\xi}_k^{\Delta} \defeq \boldsymbol{\xi}_k^{\mathrm{a}} - \boldsymbol{\xi}_k^{\mathrm{b}}$ denote the incremental dynamics for any variable $\boldsymbol{\xi}_k$.
We analyze the component-wise evolution of the state difference for any $j \in \mathcal{J}$, based on the state update equation~\eqref{eq:CFN_layer_state_update}:
\begin{equation*}
    h_{j,k+1}^{\Delta} \! = \! \left( f_{j,k}^{\mathrm{a}} \tanhfun(h_{j,k}^{\mathrm{a}}) + i_{j,k}^{\mathrm{a}} \tilde{h}_{j,k}^{\mathrm{a}} \right) \! - \! \left( f_{j,k}^{\mathrm{b}} \tanhfun(h_{j,k}^{\mathrm{b}}) + i_{j,k}^{\mathrm{b}} \tilde{h}_{j,k}^{\mathrm{b}} \right)\!.
\end{equation*}
Adding and subtracting the cross-terms $f_{j,k}^{\mathrm{a}} \tanhfun(h_{j,k}^{\mathrm{b}})$ and $i_{j,k}^{\mathrm{a}} \tilde{h}_{j,k}^{\mathrm{b}}$, we can directly group the terms as follows:
\begin{tequation}
\label{eq:diff_grouped}
    \begin{taligned}
    h_{j,k+1}^{\Delta} \! &= \! \underbrace{ f_{j,k}^{\mathrm{a}} \! \left( \tanhfun(h_{j,k}^{\mathrm{a}}) \! -  \! \tanhfun(h_{j,k}^{\mathrm{b}}) \right) }_{\tau_1} \! + \! \underbrace{ \left( f_{j,k}^{\mathrm{a}} \! - \! f_{j,k}^{\mathrm{b}} \right) \! \tanhfun(h_{j,k}^{\mathrm{b}}) }_{\tau_2} \vphantom{\underbrace{\dots}_{\tau_3}}+\\
    &\quad + \underbrace{ i_{j,k}^{\mathrm{a}} \left( \tilde{h}_{j,k}^{\mathrm{a}} - \tilde{h}_{j,k}^{\mathrm{b}} \right) }_{\tau_3} + \underbrace{ \left( i_{j,k}^{\mathrm{a}} - i_{j,k}^{\mathrm{b}} \right) \tilde{h}_{j,k}^{\mathrm{b}} }_{\tau_4}.
    \end{taligned}
\end{tequation}
We bound each term individually leveraging $|i_{j,k}^{\mathrm{a}}| \leq 1$, $|\tanhfun(h_{j,k}^{\mathrm{b}})| \leq 1$, $|f_{j,k}^{\mathrm{a}}| \leq \bar{\sigma}_f$ from~\eqref{eq:sigma_bar_deriv}, and $|\tilde{h}_{j,k}^{\mathrm{b}}| \leq \bar{\tanhfun}_{\tilde{h}}$ from~\eqref{eq:tanh_bar_deriv}.
Applying the 1-Lipschitz continuity property of the hyperbolic tangent function and the $\frac{1}{4}$-Lipschitz continuity property of the sigmoid function to the gates, we bound the absolute values of the incremental components:
\begin{equation*} 
\label{eq:tau_bound}
    \begin{taligned}
   \abs{\tau_1} &\leq \bar{\sigma}_f \infnorm{\boldsymbol{h}_k^{\Delta}}, \ \abs{\tau_3} \leq \infnorm{W_{\tilde{h}}}\infnorm{\boldsymbol{\tilde{u}}_k^{\Delta}}, \\
   \abs{\tau_2} &\leq \frac{1}{4} \left( \infnorm{W_f}\infnorm{\boldsymbol{\tilde{u}}_k^{\Delta}} + \infnorm{R_f}\infnorm{\boldsymbol{h}_k^{\Delta}} \right), \\
   \abs{\tau_4} &\leq \frac{1}{4} \left( \infnorm{W_i}\infnorm{\boldsymbol{\tilde{u}}_k^{\Delta}} + \infnorm{R_i}\infnorm{\boldsymbol{h}_k^{\Delta}} \right) \bar{\tanhfun}_{\tilde{h}}.
    \end{taligned} 
\end{equation*}
Substituting these bounds into~\eqref{eq:diff_grouped} and taking the maximum with respect to $j \in \mathcal{J}$, we obtain:
\begin{tequation}
\label{eq:diss_contraction}
    \infnorm{\boldsymbol{h}_{k+1}^{\Delta}} \leq \rho \infnorm{\boldsymbol{h}_k^{\Delta}} + \Gamma_{\tilde{u}} \infnorm{\boldsymbol{\tilde{u}}_k^{\Delta}},
\end{tequation}
where $\rho \defeq \bar{\sigma}_{f} + \frac{1}{4} \infnorm{R_{f}} + \frac{1}{4} \infnorm{R_{i}} \bar{\tanhfun}_{\tilde{h}}$ matches the left-hand side of~\eqref{eq:CFN_layer_state_update_dISS_condition}, and $\Gamma_{\tilde{u}} \defeq \infnorm{W_{\tilde{h}}} + \frac{1}{4}\infnorm{W_f} + \frac{1}{4}\infnorm{W_i}\bar{\tanhfun}_{\tilde{h}}$.
Assuming Condition~\eqref{eq:CFN_layer_state_update_dISS_condition} holds ($\rho < 1$) and unwinding the recursion back to the initial state difference, we obtain:
\begin{equation*}
    \begin{taligned}
    \infnorm{\boldsymbol{h}_{k}^{\Delta}} &\leq \rho^k \infnorm{\boldsymbol{h}_{0}^{\Delta}} + \Gamma_{\tilde{u}} \sum_{z=0}^{k-1} \rho^{k-1-z} \infnorm{\boldsymbol{\tilde{u}}_z^{\Delta}} \\
    &\leq \rho^k \infnorm{\boldsymbol{h}_{0}^{\Delta}} + \frac{\Gamma_{\tilde{u}}}{1-\rho} \max_{0 \leq z < k} \infnorm{\boldsymbol{\tilde{u}}_z^{\Delta}},
    \end{taligned}
\end{equation*}
where the input sequence is bounded by its maximum and the geometric series is bounded as in the proof of Theorem~\ref{th:CFN_layer_ISS}.
We can now deduce the $\mathcal{KL}$ and the $\mathcal{K}_\infty$ functions for the \dISS{} condition in~\eqref{eq:diss}, i.e.,
\begin{equation*}
    \begin{taligned}
    \beta_\delta\left(\infnorm{\boldsymbol{h}_{0}^{\mathrm{a}} - \boldsymbol{h}_{0}^{\mathrm{b}}}, k\right) &\defeq \rho^k \infnorm{\boldsymbol{h}_{0}^{\mathrm{a}} - \boldsymbol{h}_{0}^{\mathrm{b}}}, \\
    \gamma_{\delta u}\left(\max_{0 \leq z < k} \infnorm{\boldsymbol{\tilde{u}}^{\mathrm{a}}_z - \boldsymbol{\tilde{u}}^{\mathrm{b}}_z}\right) &\defeq \frac{\Gamma_{\tilde{u}}}{1-\rho} \max_{0 \leq z < k} \infnorm{\boldsymbol{\tilde{u}}_z^{\mathrm{a}} - \boldsymbol{\tilde{u}}_z^{\mathrm{b}}}.
    \end{taligned}
\end{equation*}
Thus, the $l$-th CFN layer, $l \in \layers$, is \dISS{} in $\mathcal{H}_{\mathrm{inv}}^{(l)}$ with respect to $\tilde{\mathcal{U}}^{(l)}$ if Condition~\eqref{eq:CFN_layer_state_update_dISS_condition} holds.
\hfill$\blacksquare$
%

\paragraphfont{Proof of Theorem~\ref{th:network_ISS}}
Every CFN layer is unconditionally ISS (Theorem~\ref{th:CFN_layer_ISS}).
The entire network is a cascade of ISS subsystems.
As proven in~\cite{jiang_input--state_2001}, the cascade of ISS systems is ISS.
Thus, the network is ISS.
\hfill$\blacksquare$
%

\paragraphfont{Proof of Theorem~\ref{th:network_dISS}}
Consider two state trajectories $\boldsymbol{h}_k^{\mathrm{a}}, \boldsymbol{h}_k^{\mathrm{b}}$ for the $L$-layer CFN driven by the input sequences $\{\boldsymbol{u}_z^{\mathrm{a}} \in \mathcal{U}\}_{z=0}^{k-1}$ and $\{\boldsymbol{u}_z^{\mathrm{b}} \in \mathcal{U}\}_{z=0}^{k-1}$, starting from initial states $\boldsymbol{h}_0^{\mathrm{a}}, \boldsymbol{h}_0^{\mathrm{b}} \in \mathcal{H}_{\mathrm{inv}}$.
As in the proof of Theorem~\ref{th:CFN_layer_dISS}, the superscript $\Delta$ denotes the incremental dynamics.
Let $\boldsymbol{\eta}_{k}^{\Delta} \in \realset_{\geq0}^{L}$ be the vector collecting the state difference norms, i.e., $\boldsymbol{\eta}_k^{\Delta} \defeq{} \left[\infnorm{\boldsymbol{h}_{k}^{(1), \Delta}}, \ldots, \infnorm{\boldsymbol{h}_{k}^{(L), \Delta}}\right]^\top$.
Recall the inter-layer connections $\infnorm{\boldsymbol{\tilde{u}}_{k}^{(l),\Delta}} = \infnorm{\boldsymbol{h}_{k+1}^{(l-1),\Delta}}$ for $l \in \layers \setminus \{1\}$ in~\eqref{eq:layer_input_def}, while $\infnorm{\boldsymbol{\tilde{u}}_{k}^{(1),\Delta}} = \infnorm{\boldsymbol{u}_{k}^{\Delta}}$.
Recursively unwrapping~\eqref{eq:diss_contraction} through the layers~\cite[Appendix B]{bonassi_stability_2021} yields the overall vector dynamics ($\leq$ is applied element-wise):
\begin{tequation}
    \label{eq:network_vector_dynamics}
    \boldsymbol{\eta}_{k}^{\Delta} \leq A_\delta \boldsymbol{\eta}_{k-1}^{\Delta} + B_{\delta u} \infnorm{\boldsymbol{u}_{k-1}^{\Delta}},
\end{tequation}
where $A_\delta \in \realset_{\geq 0}^{L \times L}$ is lower triangular with diagonal elements $\rho^{(l)}$, defined in~\eqref{eq:CFN_layer_state_update_dISS_condition}, and $B_{\delta u} \in \realset_{\geq 0}^{L}$.
%
If Condition \eqref{eq:CFN_layer_state_update_dISS_condition} holds $\forall l \in \layers$, $A_\delta$ is Schur stable as its eigenvalues (its diagonal elements) satisfy $|\rho^{(l)}| < 1$.
Propagating~\eqref{eq:network_vector_dynamics} backwards to $k=0$ yields:
\begin{equation*}
    \boldsymbol{\eta}_{k}^{\Delta} \leq A_\delta^k \boldsymbol{\eta}_{0}^{\Delta} + \sum_{z = 0}^{k - 1} A_{\delta}^{k - 1 - z} B_{\delta u} \max_{0 \leq z < k} \infnorm{\boldsymbol{u}_{z}^{\Delta}}.
\end{equation*}
Since $A_{\delta}$ and $B_{\delta u}$ are non-negative, the finite sum is upper-bounded by the infinite series. Further, due to Schur stability, we have $\sum_{z=0}^{\infty} A_{\delta}^z = (I - A_\delta)^{-1}$~\cite[Ch. 7]{meyer2023matrix}. 
Then, $\sum_{z = 0}^{k - 1} A_{\delta}^{k - 1 - z} B_{\delta u} \max_{0 \leq z < k} \infnorm{\boldsymbol{u}_{z}^{\Delta}} \leq \left(I - A_{\delta} \right)^{-1} B_{\delta u} \max_{0 \leq z < k} \infnorm{\boldsymbol{u}_{z}^{\Delta}}$.
Applying the infinity norm to both sides and noting $\infnorm{\boldsymbol{\eta}_k^{\Delta}} = \infnorm{\boldsymbol{h}_k^{\Delta}}$ results in:
\begin{equation*}
    \infnorm{\boldsymbol{h}_k^{\Delta}} \leq \infnorm{A_{\delta}^k} \infnorm{\boldsymbol{h}_0^{\Delta}} + \infnorm{\left(I - A_{\delta}\right)^{-1} B_{\delta u}} \max_{0 \leq z < k} \infnorm{\boldsymbol{u}_{z}^{\Delta}}.
\end{equation*}
Next, we leverage the Jordan form of $A_{\delta}$ to derive an upper bound for $\infnorm{A_\delta^k}$~\cite[Ch. 7]{meyer2023matrix}:
\begin{tequation}
\label{eq:jordan_bound}
    \infnorm{A_{\delta}^k} \leq \omega \sum_{z = 0}^{m - 1} \binom{k}{z} \left(\max_{l\in\layers}\rho^{(l)}\right)^{k-z},
\end{tequation}
for some $\omega \geq 1$, where $m \in \posintegerset$, $m \leq L$, is the size of the largest Jordan block of $A_{\delta}$. 
Conventionally, we set $\binom{k}{z} = 0$ for any $z > k$. 
Since $\rho^{(l)} \in (0, 1)$, the bound in~\eqref{eq:jordan_bound} tends to $0$ as $k \to \infty$.
%
Thus, the \dISS{} condition in~\eqref{eq:diss} holds with:
\begin{talign*}
    &\beta_{\delta} \! \left(\|\boldsymbol{h}_{0}^{\mathrm{a}} - \boldsymbol{h}_{0}^{\mathrm{b}}\|, k\right) \! \! \defeq \! \omega \sum_{z = 0}^{m - 1} \! \! \binom{k}{z} \! \! \left(\max_{l\in\layers}\rho^{(l)}\right)^{\!\!k-z}  \!\!\! \|\boldsymbol{h}_{0}^{\mathrm{a}} - \boldsymbol{h}_{0}^{\mathrm{b}}\|, \\
    &\gamma_{\delta u} \! \! \left( \! \max_{0 \leq z < k} \! \! \|\boldsymbol{u}_{z}^{\mathrm{a}} - \boldsymbol{u}_{z}^{\mathrm{b}}\| \! \!  \right) \! \! \defeq \! \! \infnorm{(I - A_\delta)^{-1} \! B_{\delta u}} \! \max_{0 \leq z < k}\! \!  \|\boldsymbol{u}_{z}^{\mathrm{a}} - \boldsymbol{u}_{z}^{\mathrm{b}}\|.
\end{talign*}
Hence, the $L$-layer CFN is \dISS{} in $\mathcal{H}_{\mathrm{inv}}$ with respect to $\mathcal{U}$ if Condition~\eqref{eq:CFN_layer_state_update_dISS_condition} holds for every layer $l \in \layers$.
\hfill$\blacksquare$
%

\paragraphfont{Proof of Theorem~\ref{th:DGN_stability}}
The DGN is a special case of the CFN where $R_f^{(l)} = R_i^{(l)} = \boldsymbol{0}, \forall l \in \layers$.
Substituting these zero matrices into the \dISS{} sufficient condition in~\eqref{eq:CFN_layer_state_update_dISS_condition}, the left-hand side simplifies to $\rho^{(l)} = \bar{\sigma}_{f}^{(l)}$.
%
Since $\bar{\sigma}_{f}^{(l)} \in (0, 1)$ as in~\eqref{eq:sigma_bar_deriv}, the condition $\rho^{(l)} < 1$ is satisfied.
Since this holds $\forall l \in \layers$, the $L$-layer DGN is \dISS{} by Theorem~\ref{th:network_dISS}.
\hfill$\blacksquare$
\bibliographystyle{IEEEtran}
\bibliography{bibliography}
\end{document}

%% file: figures/cfn-dgn.tex
\tikzmath{
    \lblock = 0.6 cm;
    \radius = \lblock;
    \eqwidth = 3 cm; 
    \rnnheight = \lblock;
    \grid = 0.5;       
    \halfgrid = \grid/2;
}

\tikzset{
    block/.style={
        rectangle,
        rounded corners,
        draw,
        thick,
        minimum width=\lblock,
        minimum height=\lblock,
        align=center
    },
    eqblock/.style={
        rectangle,
        rounded corners,
        draw,
        thick,
        minimum width=\eqwidth,
        minimum height=\rnnheight,
        align=center
    },
    surround/.style={
        rectangle,
        draw,
        thick
    },
    arrow/.style={
        ->,
        >=stealth,
        thick
    },
    line/.style={
        thick
    },
    optional/.style={
        dashed,
        line cap=round,
        dash pattern=on 3pt off 3pt,
        draw=state_color
    }
}

\begin{tikzpicture}[every node/.style={font=\small}]

    
    \node[block, draw=input_color, text=input_color] (uk) {$\boldsymbol{\tilde{u}}_k^{(l)}$};
    
    \node[block, draw=state_color, text=state_color, above=0.75cm of uk.north] (hk) {$\boldsymbol{h}_k^{(l)}$};

    \node[block, right=\grid cm of hk.east] (fk) {$\boldsymbol{f}_k^{(l)}$};
    \node[block, right=\grid cm of uk.east] (htk) {$\boldsymbol{\tilde{h}}_k^{(l)}$};
    
    \path (fk) -- (htk) node[midway, block] (ik) {$\boldsymbol{i}_k^{(l)}$};

    \node[eqblock, right=\halfgrid cm of ik.east] (eq) {
        $\textcolor{state_color}{\boldsymbol{h}_{k+1}^{(l)}} = \boldsymbol{f}_k^{(l)} \circ \boldsymbol{\tanhfun}\left(\textcolor{state_color}{\boldsymbol{h}_k^{(l)}}\right) + \boldsymbol{i}_k^{(l)} \circ \boldsymbol{\tilde{h}}_k^{(l)}$
    };

    \node[block, draw=state_color, text=state_color, right=\halfgrid cm of eq.east] (hkp) {$\boldsymbol{h}_{k+1}^{(l)}$};


    \draw[arrow, draw=input_color] (uk) -- (htk);
    \draw[arrow, draw=input_color] ($ (uk.east)!0.5!(htk.west) $) |- (ik.west |- {$(ik.south)!0.3!(ik.north)$});
    \draw[arrow, draw=input_color] ($ (uk.east)!0.5!(htk.west) $) |- (fk.west |- {$(fk.south)!0.3!(fk.north)$});

    
    \draw[arrow, optional] (hk) -- (fk);
    
    \draw[arrow, optional] ($ (hk.east)!0.3!(fk.west) $) |- (ik.west |- {$(ik.north)!0.5!(ik.south)$});
    
    \draw[arrow, draw=state_color] ($ (hk.east)!0.3!(fk.west) $) 
        -- ++(0, 0.4cm) 
        -| (eq.north -| {$(eq.west)!0.2!(eq.east)$});

    \draw[arrow] (fk.east) -- ++(0.2, 0) -| (eq.north -| {$(eq.west)!0.1!(eq.east)$});
    \draw[arrow] (ik) -- (eq);
    \draw[arrow] (htk.east) -- ++(0.2, 0) -| (eq.south -| {$(eq.west)!0.1!(eq.east)$});

    \draw[arrow] (eq) -- (hkp);

\end{tikzpicture}

%% file: figures/realizations.tex
\begin{tikzpicture}
    \pgfplotsset{
        real_style/.style={
            draw=black,
            very thick,
            densely dashed,
            smooth
        },
        cfn_style/.style={
            draw=cfn_color, 
            line width=2pt,
            opacity=0.5,
            solid,
            smooth,
            line join=round
        },
        dgn_style/.style={
            draw=dgn_color, 
            line width=1.0pt,
            solid,
            smooth,
            line join=round
        }
    }

    \begin{groupplot}[
        group style={
            group size=1 by 2,
            horizontal sep=0cm,
            vertical sep=0.5cm, 
            x descriptions at=edge bottom,
            y descriptions at=edge left,
        },
        width=\columnwidth,
        height=3.5cm,
        grid=major,
        each nth point={2},
        filter discard warning=false,
        unbounded coords=jump,
        legend style={
            at={(0.5,1.05)}, 
            anchor=south, 
            legend columns=3, 
            draw=none, 
            fill=none,
            font=\footnotesize
        },
        xmin=0,
    ]

    \nextgroupplot[
        ylabel={pH level [-]},
        xmax= 5.5, 
        xtick={0, 1, 2, 3, 4, 5},
        xticklabels={0, 1, 2, 3, 4, 5},
        ytick={6, 7, 8},
    ]
        
        \addplot[cfn_style] table [x=k, y=y_cfn] {data/dataset_1.dat};
        \addlegendentry{CFN}

        \addplot[dgn_style] table [x=k, y=y_dgn] {data/dataset_1.dat};
        \addlegendentry{DGN}

        \addplot[real_style] table [x=k, y=y_real] {data/dataset_1.dat};
        \addlegendentry{Real}

    \nextgroupplot[
        ylabel={Water level [m]},
        xmax= 20.8, 
        xtick={0, 5, 10, 15, 20},
        ytick={0, 0.5, 1, 1.5},
        xlabel={Time [h]},
        legend style={draw=none, fill=none, at={(0.5,1.1)}, anchor=south}, 
        legend to name=sharedlegend 
    ]
        
        \addplot[cfn_style] table [x=k, y=y_cfn] {data/dataset_2.dat};
        
        \addplot[dgn_style] table [x=k, y=y_dgn] {data/dataset_2.dat};

        \addplot[real_style] table [x=k, y=y_real] {data/dataset_2.dat};

    \end{groupplot}
\end{tikzpicture}